\documentclass{amsart}
\usepackage{amssymb}
\usepackage{amsmath}
\usepackage{mathrsfs}                  % for \mathscr
\usepackage{appendix}
\usepackage[dvips,ps,all,color]{xy}    % diagrams with color
\UseCrayolaColors

\title%
[Integrality results]
{The ideal of relations for the ring of invariants\\of $n$ points on
the line:  integrality results}

\author[Benjamin Howard, John Millson, Andrew Snowden and Ravi Vakil]
{Benjamin Howard, John Millson, Andrew Snowden and Ravi Vakil*}

\thanks{*B.\ Howard was supported by NSF fellowship DMS-0703674.
J.\ Millson was supported by the NSF grant DMS-0405606, the NSF FRG grant
DMS-0554254 and the Simons Foundation.
A. Snowden was partially supported by NSF fellowship DMS-0902661.
R.\ Vakil was partially supported by NSF grant DMS-0801196}

\date{September 17, 2009.}

\newcommand{\cut}[1]{}

\newcommand{\Ncom}[1]{}
\newcommand{\Rnts}[1]{}
\newtheorem{theorem}{Theorem}[section]
\newtheorem{lemma}[theorem]{Lemma}
\newtheorem{corollary}[theorem]{Corollary}
\newtheorem{proposition}[theorem]{Proposition}

\theoremstyle{definition}
\newtheorem{definition}[theorem]{Definition}

\theoremstyle{remark}
\newtheorem{remark}[theorem]{\it Remark}

\newcommand{\reg}{\mathrm{reg}}

\DeclareMathOperator{\Sym}{Sym}

\DeclareMathOperator{\gr}{gr}

\DeclareMathOperator{\Hom}{Hom}

\def\SL{\mathrm{SL}}

\def\Z{\mathbb{Z}}

\def\P{\mathbb{P}}

\def\Q{\mathbb{Q}}

\def\half{\tfrac{1}{2}}

\def\wt{\mathbf{wt}}
\def\cq{/\!/}

\let\ms\mathscr

\let\ss\scriptstyle
\let\wt\widetilde
\let\ol\overline

\newcommand\dedge[1]{\overrightarrow{#1\rule{0em}{.7em}}}
\newcommand\uedge[1]{\overline{#1\rule{0em}{.7em}}}

\begin{document}

    %%%%%%%%%%%%%%%%%%%%
    %% ARTICLE HEADER %%
    %%%%%%%%%%%%%%%%%%%%

%abstract
\begin{abstract}
Consider the projective coordinate ring of the GIT quotient $(\P^1)^n \cq
\SL(2)$, with the usual linearization, where $n$ is even.  In 1894, Kempe
proved that this ring is generated in degree one.  In \cite{hmsv2} we showed
that, over $\Q$, the relations between degree one invariants are generated by a 
class of quadratic relations --- the \emph{simplest binomial relations} ---
with the exception of $n=6$, where there is a single cubic relation.  The
purpose of this paper is to show that these results hold over
$\Z[\tfrac{1}{12!}]$, and to suggest why they may be true over
$\Z[\tfrac 1 6]$.
\end{abstract}

%% title and abstract
\maketitle

%% table of contents
\tableofcontents

\section{Introduction}

Let $L$ be a set of even cardinality $n$ (what we will call an {\em even set}).
Let $R_L$ denote the projective coordinate ring of the space of $L$ points on
$\P^1$, modulo the diagonal action of $\SL(2)$:
\begin{displaymath}
R_L=\bigoplus_{k=0}^{\infty} R_L^{(k)}, \qquad
R_L^{(k)}=\Gamma(X_L, \ms{O}_{X_L}(k))^{\SL(2)}
\end{displaymath}
where $X_L=\Hom(L, \P^1)$ and $\ms{O}_{X_L}(1)$ is the exterior tensor product
of $L$ copies of $\ms{O}_{\P^1}(1)$.  Kempe showed 
that $R_L$ is generated by its first graded piece $V_L=R_L^{(1)}$; a proof
can be found in \cite[\S 2.3]{hmsv2}.  Let $I_L$ be the \emph{ideal of
relations} between degree one invariants, the kernel of the surjection
$\Sym(V_L) \to R_L$.

In \cite[Thm.~9.1]{hmsv2}, we
showed that the
quadratic \emph{simplest binomial relations} (elements of $I_L^{(2)}$,
defined in \S \ref{simple}) generate $I_L$ over $\Z[\tfrac{1}{n!}]$ for $n \ge 8$.
The main result of this paper is the following:

\begin{theorem}[Main theorem]
\label{mainthm}
Let $L$ be an even set of cardinality at least 8.  Then $I_L$ is
generated over $\Z[\tfrac{1}{12!}]$ by the simplest binomial relations.
\end{theorem}

In \cite[\S 8]{hmsv2} we proved that $I_L$ is generated over $\Z[\tfrac{1}{12!}]$
by quadratics.  We also proved in \cite[\S 9]{hmsv2} that the ``simple
binomial relations'' (defined in \S \ref{simple}) span
the same space as the simplest binomial relations over $\Z[\half]$.  Thus to
prove Theorem~\ref{mainthm} it suffices to prove:

\begin{theorem}
\label{mainthm2}
Let $L$ be an even set.  Then $I_L^{(2)}$ is spanned over
$\Z[\tfrac{1}{12!}]$ by the simple binomial relations.
\end{theorem}

We established this statement over $\Z[\tfrac{1}{n!}]$ in
\cite[\S 9]{hmsv2} by a
short argument using the representation theory of the symmetric
group $S_n$.  The use of representation theory is the reason that we required
$n!$ to be invertible.  In this paper, we prove the result over
$\Z[\tfrac{1}{12!}]$ by a totally different --- and much more involved ---
inductive argument.  The work lies in establishing the inductive step:

\begin{theorem}
\label{mainthm3}
Let $L$ be an even  set of cardinality at least 10 but not 12.  Assume
that for all sets $L'$ of smaller even cardinality $I_{L'}^{(2)}$ is generated
over $\Z[\tfrac{1}{2N}]$ by simple binomial relations.  Then the same is
true for $I_L^{(2)}$.
\end{theorem}

As remarked above, the results of \cite{hmsv2} show that $I_L^{(2)}$ is spanned
by simple binomial relations over $\Z[\tfrac{1}{12!}]$ whenever $L$ has
cardinality at most 12.  Thus Theorem~\ref{mainthm3} shows that this statement
holds for $L$ of any cardinality, which establishes Theorem~\ref{mainthm2}.

\begin{remark}
The reader may wonder about the importance of the simplest binomial
relations.  Aside from aesthetic reasons, 
there is the following one:  the simplest binomial relations form a single
orbit under the action of the symmetric group.  Thus Theorem~\ref{mainthm}
shows that $I_L$ is a ``principal $S_n$-ideal'' over
$\Z[\tfrac{1}{12!}]$.  Perhaps most important, they naturally arise in
much 
more general circumstances in the retrogeneration conjecture of the
third author.
\end{remark}

\subsection{Improving Theorem~\ref{mainthm}}

We suspect  that Theorem~\ref{mainthm} may in fact hold over $\Z[\tfrac{1}{6}]$.
Here are some comments along these lines:
\begin{itemize}
\item When $L$ has cardinality 8, Theorem~\ref{mainthm} is true over
$\Z[\tfrac{1}{6}]$.  In fact, $I_L^{(2)}$ is then spanned over
$\Z[\half]$ by the simplest binomial relations (this follows from
\cite[Proposition~2.10]{hmsv1} which gives a computer aided proof that the
simple binomial relations span $I^{(2)}_L$ over $\Z$) while $I_L$ is
generated over $\Z[\tfrac{1}{3}]$ by quadratics (see \cite{8pts}).
\item Let $L$ have cardinality 10.  Then Theorem~\ref{mainthm2} shows that the
simplest binomial relations span $I_L^{(2)}$ over $\Z[\tfrac{1}{2}]$, since
this is the case for all smaller cardinalities.  However, we cannot conclude
from this that $I_L^{(2)}$ is generated by simplest binomial
relations over $\Z[\tfrac{1}{6}]$ because we only know that $I_L$ is generated
by quadratics over $\Z[\frac{1}{10!}]$.
\item By the results of \cite{hmsv2}, to prove Theorem~\ref{mainthm} over
$\Z[\tfrac{1}{6}]$ it would suffice to verify the following two statements:
(1) $I_L^{(3)}$ is generated by quadratics (i.e., the map $V_L \otimes
I_L^{(2)} \to I_L^{(3)}$ is surjective) over $\Z[\tfrac{1}{6}]$ when $L$ has
cardinality 10 or 12; (2) $I_L^{(2)}$ is spanned over $\Z[\tfrac{1}{6}]$ by
simple binomial relations when $L$ has cardinality 12.
\end{itemize}
The above comments show that to confirm our suspicion that
Theorem~\ref{mainthm} holds over $\Z[\tfrac{1}{6}]$ one needs only 
perform a few calculations.  However, these computations are quite large: for
instance, when $L$ has cardinality 12, the space $I_L^{(3)}$ has
dimension  339,240.   Thus to carry
this out on a computer would require 
considerable thought.

\subsection*{Acknowledgments}

We thank Chris Manon and Lawrence O'Neil for some helpful
discussions.

\section{Background}

\subsection{The ring and the ideal}

We now recall the graphical description of the ring $R_L$.  For more details
see \cite[\S 2]{hmsv2}.  Let $\ms{G}_L^{\reg}$ denote the set of all regular
directed graphs on the vertex set $L$.  ``Regular'' means that each vertex
has the same valence.  We allow our graphs to have loops (edges from a vertex
to itself) and multiple edges between the same vertices.  We give
$\ms{G}_L^{\reg}$ the structure of a semi-group by defining $\Gamma \cdot
\Gamma'$ to be the graph on $L$ whose edge set is the disjoint union of those
of $\Gamma$ and $\Gamma'$.  For a graph $\Gamma$ we let $X_{\Gamma}$ denote
the corresponding element of the semi-group algebra $\Z[\ms{G}^{\reg}_L]$.

The ring $R_L$ is defined to be the quotient of $\Z[\ms{G}^{\reg}_L]$ by the
following three types of relations:
\begin{itemize}
\item \emph{Loop relation:} $X_{\Gamma}=0$ if $\Gamma$ has a loop.
\item \emph{Sign relation:} $X_{\Gamma}=-X_{\Gamma'}$ if $\Gamma'$ is obtained
from $\Gamma$ be reversing the direction of a single edge.
\item \emph{Pl\"ucker relation:} Let $\Gamma$ be an element of
$\ms{G}_L^{\reg}$ and let $\dedge{ab}$ and $\dedge{cd}$ be two edges of
$\Gamma$.  Let $\Gamma'$ (resp.\ $\Gamma''$) be the graph obtained by
replacing $(\dedge{ab}, \dedge{cd})$ in $\Gamma$ with $(\dedge{ad},
\dedge{cd})$ (resp.\ $(\dedge{ac}, \dedge{bd})$).  Then $X_{\Gamma}=X_{\Gamma'}
+X_{\Gamma''}$.
\end{itemize}
We denote the image of $X_{\Gamma}$ in $R_L$ still by $X_{\Gamma}$.  We
give $R_L$ a grading by letting $R_L^{(k)}$ be the span of the $X_{\Gamma}$
with $\Gamma$ of regular degree $k$.

Assume now that $L$ is an even set.  The ring $R_L$ is then generated by its first graded piece
$V_L=R_L^{(1)}$ (Kempe's theorem).  The space $V_L$ is spanned by those
$X_{\Gamma}$ where $\Gamma$ is a \emph{matching}, that is, a regular graph of
degree one.  We let $I_L$ be the \emph{ideal of relations}, that is, the kernel
of the surjection $\Sym(V_L) \to R_L$.  We let $W_L$ be the second graded
piece of $R_L$, namely $R_L^{(2)}$, and we let $B_L$ be the kernel of the
map $V_L^{\otimes 2} \to W_L$.  The sequence
\begin{displaymath}
0 \to B_L \to V_L^{\otimes 2} \to W_L \to 0
\end{displaymath}
is exact.  The space of quadratic relations $I_L^{(2)}$ is the image of $B_L$ under the map
$V_L^{\otimes 2} \to \Sym^2(V_L)$.

\subsection{Simple binomial relations}

\label{simple}Let $U$ be a subset of $L$ of cardinality four and put $L'=L \setminus U$.
Let $\Delta$ and $\Delta'$ be matchings on $U$ and let $\Gamma$ and $\Gamma'$
be matchings on $L'$.  We then have the following obvious element of
$I_L^{(2)}$:
\begin{displaymath}
X_{\Gamma \Delta} X_{\Gamma' \Delta'} - X_{\Gamma \Delta'} X_{\Gamma' \Delta}.
\end{displaymath}
(Here $\Gamma \Delta$, etc., denotes the union of $\Gamma$ and $\Delta$, a
matching with vertex set $L$.)  We regard this as an element of $B_L$ as
well by inserting the $\otimes$ symbol in between the $X$'s in the above.
We call these elements the \emph{simple binomial relations}.  As remarked in
the introduction, these were shown in \cite{hmsv2} to span $I_L^{(2)}$ when
certain primes are inverted; the goal of this paper is to show that one
needs to invert fewer primes than was done there.  (The ``simplest binomial
relations'' referred to in the introduction are the simple binomial relations
where $\Gamma$ is a union of a single 4-cycle with a bunch of 2-cycles; they
will not play an important role in the paper.)

Pure tensors in $V_L^{\otimes 2}$ can be thought of as graphs on $L$ in which
each edge has been assigned one of two colors in such a manner that each
vertex belongs to exactly one edge of each color.  Elements of $\Sym^2(V_L)$
can be thought of in a similar manner.  See \cite[\S 5.1]{hmsv2} for a more
complete discussion.  With this mode of thinking, one can depict the simple
binomial relations in a very visually pleasing way.  For example, here is one
on 10 points:
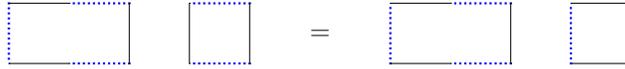
\begin{displaymath}
\begin{xy}
(0, 4)*{}="A"; (8, 4)*{}="B"; (16, 4)*{}="C"; (24, 4)*{}="D";
(0, -4)*{}="E"; (8, -4)*{}="F"; (16, -4)*{}="G"; (24, -4)*{}="H";
(-8, 4)*{}="P"; (-8, -4)*{}="Q";
{\ar@[blue]@{..}@[|(2)] "A"; "B"};
{\ar@[blue]@{..}@[|(2)] "Q"; "P"};
{\ar@[blue]@{..}@[|(2)] "E"; "F"};
{\ar@[black]@{-} "B"; "F"};
{\ar@[black]@{-} "P"; "A"};
{\ar@[black]@{-} "Q"; "E"};
{\ar@[blue]@{..}@[|(2)] "C"; "D"};
{\ar@[blue]@{..}@[|(2)] "G"; "H"};
{\ar@[black]@{-} "C"; "G"};
{\ar@[black]@{-} "D"; "H"};
\end{xy}
\qquad = \qquad
\begin{xy}
(0, 4)*{}="A"; (8, 4)*{}="B"; (16, 4)*{}="C"; (24, 4)*{}="D";
(0, -4)*{}="E"; (8, -4)*{}="F"; (16, -4)*{}="G"; (24, -4)*{}="H";
(-8, 4)*{}="P"; (-8, -4)*{}="Q";
{\ar@[blue]@{..}@[|(2)] "A"; "B"};
{\ar@[blue]@{..}@[|(2)] "E"; "F"};
{\ar@[blue]@{..}@[|(2)] "P"; "Q"};
{\ar@[black]@{-} "P"; "A"};
{\ar@[black]@{-} "Q"; "E"};
{\ar@[black]@{-} "B"; "F"};
{\ar@[black]@{-} "C"; "D"};
{\ar@[black]@{-} "G"; "H"};
{\ar@[blue]@{..}@[|(2)] "C"; "G"};
{\ar@[blue]@{..}@[|(2)] "D"; "H"};
\end{xy}
\end{displaymath}
The 4-cycles are formed out of $\Delta$ and $\Delta'$ while the 6-cycles are
formed out of $\Gamma$ and $\Gamma'$.

There is a natural map $\boxtimes:V_L \otimes V_{L'} \to V_{L \amalg L'}$
taking $X_{\Gamma} \otimes X_{\Delta}$ to $X_{\Gamma \Delta}$.  We called
this \emph{outer multiplication} in \cite[\S 7.1]{hmsv2}.  The simple binomial
relations may be defined in terms of outer multiplication.
Let $U$ be a set of size four.  Then $B_U$ is spanned by the commutativity
relations $X_{\Delta} \otimes X_{\Delta'}- X_{\Delta'} \otimes X_{\Delta}$.
The simple binomial relations in $B_L$ are then exactly the outer products
of arbitrary elements of $V_{L'}^{\otimes 2}$ with elements of $B_U$, where
$L=L' \amalg U$.

\subsection{The straightening algorithm}

Fix an embedding of $L$ into the unit circle in the plane.  We call a
graph $\Gamma$ on $L$ \emph{planar} if its edge do not cross when drawn
as chords in the circle.  The following theorem, also due to Kempe, will be
key.  It is well-known, so we omit a proof.  For details see
\cite[Props.~2.5, 2.6]{hmsv1}.

\begin{proposition}[Kempe]
\label{prop:straight}
The $X_{\Gamma}$, with $\Gamma$ planar, span $R_L$ over $\Z$.  The only
relations between these elements are the sign relations.  Thus if for each
undirected planar graph one chooses a direction on the edges then the
resulting $X_{\Gamma}$'s form a basis for $R_L$ over $\Z$.
\end{proposition}

To write an arbitrary $X_{\Gamma}$ in terms of the planar basis one uses
the \emph{straightening algorithm}:  simply take any two edges of $\Gamma$
which cross and apply the Pl\"ucker relation to them; repeat.  This
algorithm terminates because the total lengths of the edges in the resulting
graphs is less than the total lengths of the edges in $\Gamma$, since the
lengths of the diagonals in a quadrilateral exceeds that of any pair of
opposite edges.  This shows that the planar graphs span over $\Z$.

\section{Reduction of Theorem~\ref{mainthm3} to Proposition~\ref{prop:merge}}

The goal of this section is to reduce the proof of Theorem~\ref{mainthm3} to
that of a different statement, Proposition~\ref{prop:merge}.  The statement
of this proposition is  less transparent than that of
Theorem~\ref{mainthm3}; however, in the next section we will reduce the proof
of Proposition~\ref{prop:merge} to that of a very concrete statement.

\subsection{Partitions}

Let $L$ be an even set.  By a \emph{partition} of $L$ we mean a collection
$\ms{U}$ of subsets of $L$ such that $\ms{U}$ has cardinality at least 2, each
element of $\ms{U}$ has even cardinality and $L$ is the disjoint union of the
elements of $\ms{U}$.  For a partition $\ms{U}$ we define
\begin{displaymath}
V_{\ms{U}} = \bigotimes_{U \in \ms{U}} V_U, \qquad
W_{\ms{U}} = \bigotimes_{U \in \ms{U}} W_U, \qquad
B_{\ms{U}} = \sum_{U \in \ms{U}} B_U \otimes V_{\ms{U} \setminus
\{ U \}}^{\otimes 2},
\end{displaymath}
the sum taking place inside of $V_{\ms{U}}^{\otimes 2}$ in the last
definition.  Note that the sequence
\begin{displaymath}
0 \to B_{\ms{U}} \to V_{\ms{U}}^{\otimes 2} \to W_{\ms{U}} \to 0
\end{displaymath}
is exact.  We now put
\begin{displaymath}
\wt{V}^{(2)}_L = \bigoplus_{\ms{U}} V_{\ms{U}}^{\otimes 2}, \qquad
\wt{W}_L = \bigoplus_{\ms{U}} W_{\ms{U}}, \qquad
\wt{B}_L = \bigoplus_{\ms{U}} B_{\ms{U}},
\end{displaymath}
the sum taken over all partitions $\ms{U}$.  Note that the sequence
\begin{displaymath}
0 \to \wt{B}_L \to \wt{V}^{(2)}_L \to \wt{W}_L \to 0
\end{displaymath}
is exact.

A subset $U \subset L$ is \emph{closed} with respect to a
graph $\Gamma$ if every edge of $\Gamma$ which contains a vertex in $U$ is
completely contained in $U$.  A partition $\ms{U}$ is
\emph{closed} with respect to $\Gamma$ if each of its pieces is.
One may interpret $\wt{W}_L$ as the free module on tuples $(\Gamma, \ms{U})$
where $\Gamma$ is a degree two graph on $L$ and $\ms{U}$ is a closed partition
of $L$.  One may perform Pl\"ucker relations on two edges only if
they both lie within  the same piece of the partition.
Similarly, one may interpret $\wt{V}^{(2)}_L$ as the free module on tuples $(\Gamma,
\ms{U})$ where $\Gamma$ is a colored graph of degree two such and $\ms{U}$ is
closed with respect to $\Gamma$.  One is allowed to perform Pl\"ucker
relations on two edges of the same color which lie in the same piece of the
partition.  The map $\wt{V}^{(2)}_L \to \wt{W}_L$ is then given by forgetting
the coloring.

\subsection{The main diagram}

We have maps $\wt{V}^{(2)}_L \to \wt{V}^{\otimes 2}_L$ and $\wt{W}_L \to W_L$,
given by forgetting the partition.  A key point is the following:

\begin{lemma}
Assume the cardinality of $L$ is at least six.  Then the maps $\wt{V}_L^{(2)}
\to V_L^{\otimes 2}$ and $\wt{W}_L \to W_L$ are surjective over $\Z[\half]$
\end{lemma}

\begin{proof}
It is shown in \cite[\S 6.3]{hmsv2} that $V_L^{\otimes 2}$ and $W_L$ are spanned
over $\Z[\half]$ by graphs which are unions of 2- and 4-cycles.  This can also
be seen (at least for $W_L)$ from the identities (I2) and (I3) in the appendix.
Thus if $L$ has cardinality at least 6 then one can write any element of
$V_L^{\otimes 2}$ as a sum of $X_{\Gamma_i}$ where each $\Gamma_i$ is a
2-colored graph which is disconnected, and thus belongs to some $V_{\ms{U}}$.
The same reasoning applies to $W_L$.
\end{proof}

We denote the kernel of $\wt{V}_L^{(2)} \to V_L^{\otimes 2}$ by $P_L$
and the kernel of $\wt{W}_L \to W_L$ by $Q_L$.  We now have the commutative
 diagram:
\begin{displaymath}
\xymatrix{
&& 0 \ar[d] & 0 \ar[d] \\
&& P_L \ar[r]^{\alpha_L} \ar[d] & Q_L \ar[d] \\
0 \ar[r] & \wt{B}_L \ar[d]^{\beta_L} \ar[r] & \wt{V}^{(2)}_L \ar[r]
\ar[d] & \wt{W}_L \ar[r] \ar[d] & 0 \\
0 \ar[r] & B_L \ar[r] & V_L^{\otimes 2} \ar[r] \ar[d] &
W_L \ar[r] \ar[d] & 0 \\
&& 0 & 0
}
\end{displaymath}
in which the rows are exact.  The columns are exact over $\Z[\half]$ if $L$
has cardinality at least six, by the previous lemma.

\begin{lemma}
Let $L$ be an even set of cardinality at least six.  Assume (1) $\beta_L$ is
surjective; and (2) $B_{L'}$ is spanned over $\Z[\frac{1}{2N}]$ by simple
binomial relations for all proper even subsets $L'$ of $L$.  Then $B_L$ is
spanned over $\Z[\tfrac{1}{2N}]$ by simple binomial relations.
\end{lemma}

\begin{proof}
The image of $\beta_L$ is exactly the space spanned by outer products of
elements of $B_{L'}$ and elements of $V_{L \setminus L'}^{\otimes 2}$, where
$L'$ is a proper even subset of $L$.  Hypothesis (2) ensures that $B_{L'}$ is
spanned by outer products of the commutativity relation on four points.
Thus the image of $\beta_L$ is spanned by the outer products of commutativity
relations on four points.  Since $\beta_L$ is surjective by (1), we see that
$B_L$ is spanned by outer products of commutativity relations on four points,
i.e., by simple binomial relations.
\end{proof}

\begin{lemma}
Let $L$ be an even set of cardinality at least six.  Then $\beta_L$ is
surjective over $\Z[\tfrac{1}{N}]$ if and only if $\alpha_L$ is surjective
over $\Z[\tfrac{1}{N}]$.
\end{lemma}

\begin{proof}
This is a simple diagram chase.
\end{proof}

We thus see that to prove Theorem~\ref{mainthm3} we only need to show that
$\alpha_L$ is surjective when $L$ has cardinality at least 10 but not 12.  We
will prove this by giving explicit (and simple) generators for $Q_L$ and then
verifying that these generators are in the image of $\alpha_L$.

\subsection{Generators for $Q_L$}

There is an obvious class of relations in $Q_L$, which we call the
\emph{merging relations}: the relations $(\Gamma, \ms{U})=(\Gamma, \ms{U}')$
where $\Gamma$ is a graph, $\ms{U}$ is a partition with at least three pieces
and closed with respect to $\Gamma$ and $\ms{U}'$ is obtained from $\ms{U}$ by
merging (taking the union of) two of its pieces.  We then have the following
proposition, the proof of which accounts for most of our effort.

\begin{proposition}
\label{prop:merge}
If the cardinality of $L$ is at least 10 but not 12 then the merging relations
span $Q_L$ over $\Z[\half]$.
\end{proposition}

The following lemma shows that Proposition~\ref{prop:merge} implies
Theorem~\ref{mainthm3}.

\begin{lemma}
The merging relations are in the image of the map $\alpha_L$ (over $\Z$).
\end{lemma}

\begin{proof}
Let $(\Gamma, \ms{U})-(\Gamma, \ms{U}')$ be a merging relation.  By using the
same Pl\"ucker relations on each copy of $\Gamma$, this relation can be written
as a sum of merging relations in which $\Gamma$ only has even cycles.  We can
then 2-color $\Gamma$ to get a graph $\Gamma'$ in $V_L^{\otimes 2}$.  Since we
did not change the underlying graph, $\ms{U}$ and $\ms{U}'$ are still closed
with respect to $\Gamma'$.  Thus $(\Gamma', \ms{U})-(\Gamma', \ms{U}')$ is an
element of $P_L$.  We have thus lifted the original merging relation through
$\alpha_L$.
\end{proof}

\begin{remark}
We will see below that when the cardinality of $L$ is equal to 12 there is a
single relation (up to symmetry) that one can include, called the odd cycle
exchange relation (\S \ref{oddcycle}), which together with the merging relations spans all of
$Q_L$.  However, it is not clear if this relation is in the image of
$\alpha_L$.
\end{remark}

\section{Reduction of Proposition~\ref{prop:merge} to
Proposition~\ref{prop:qp}}
\label{sec-mainthm2}

The goal of this section is to reduce the proof of Proposition~\ref{prop:merge}
to that of two much more concrete statements, given in
Proposition~\ref{prop:qp}.

\subsection{The odd cycle exchange relation}

\label{oddcycle}We first discuss a class of relations called the \emph{odd cycle exchange}
relations.  Let $\Gamma$ be a regular degree 2 graph on $L$ and let
$\{U_i\}_{i=1}^4$ be four sets of odd cardinality which are closed with
respect to $\Gamma$.  Put $\ms{U}=\{ U_1 \cup U_2, U_3 \cup U_4 \}$ and
$\ms{U}'=\{ U_1 \cup U_3, U_2 \cup U_4 \}$.  The odd cycle exchange relation
is then
\begin{displaymath}
(\Gamma, \ms{U})=(\Gamma, \ms{U}').
\end{displaymath}
Note that if the cardinality of $L$ is less than 12 then there are no
odd cycle exchange relations because $L$ cannot be split up into four sets
off odd cardinality at least 3.  On the other hand, when $L$ has cardinality
at least 14 the odd cycle exchange relations are not new:

\begin{proposition}
\label{prop:oddmerge}
If the cardinality of $L$ is at least 14 then the odd cycle exchange relations
are linear combinations of merging relations over $\Z[\half]$.
\end{proposition}

\begin{proof}
Let $\Gamma$ and $\{U_i\}_{i=1}^4$ as above be given.  First consider the
case where $\Gamma$ contains more than one cycle in one of the $U_i$, say in
$U_1$.  We can thus write $U_1=V \cup U_1'$ where $V$ has even cardinality,
$U_1'$ has odd cardinality and $U$ and $V$ are closed with respect to
$\Gamma$.  We now have the sequence of merging relations
\begin{displaymath}
\begin{split}
(\Gamma, \{ U_1 \cup U_2, \;\; U_3 \cup U_4 \})
=&(\Gamma, \{ V, \;\; U_1' \cup U_2, \;\; U_3 \cup U_4 \}) \\
=&(\Gamma, \{ V, \;\; U_1' \cup U_2 \cup U_3 \cup U_4 \}) \\
=&(\Gamma, \{V, \;\; U_1' \cup U_3, \;\; U_2 \cup U_4 \}) \\
=&(\Gamma, \{U_1 \cup U_3, \;\; U_2 \cup U_4 \})
\end{split}
\end{displaymath}
which realizes the odd cycle exchange relation.

We now handle the general case.  Since the cardinality of $L$ is at least 14,
one of the $U_i$, say $U_1$, has at least five vertices.  Now, by using only
Pl\"ucker relations in $U_1$ one may write $\Gamma=\sum a_i \Gamma_i$ where
each $\Gamma_i$ has more than one cycle in $U_1$ and $a_i$ belong to
$\Z[\half]$.  This is essentially proved in \cite[\S 6.3]{hmsv2},
but also follows from identity (I2) of the appendix.
Since we only used Pl\"ucker relations within $U_1$, all the Pl\"ucker
relations are allowable with respect to $\ms{U}=\{U_1 \cup U_2, U_3 \cup
U_4\}$ and $\ms{U}'=\{U_1 \cup U_3, U_2 \cup U_4\}$.  We thus have the
following equality in $\wt{W}_L$:
\begin{displaymath}
(\Gamma, \ms{U})-(\Gamma, \ms{U}')=\sum \left[ (\Gamma_i, \ms{U})-
(\Gamma_i, \ms{U}') \right]
\end{displaymath}
This expresses the odd cycle exchange relation in which we are interested (the
one on the left) in terms of odd cycle exchange relations in which the graph on
$U_1$ has more than one cycle.  By the previous paragraph, the right side lies
in the span of the merging relations.
\end{proof}

When the cardinality of $L$ is equal to 12 there is only one odd
cycle exchange relation, up to symmetry:
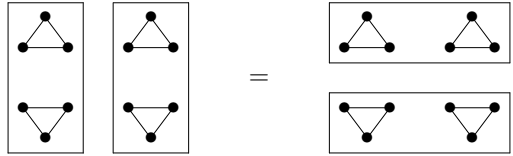
\begin{displaymath}
\begin{xy}
(-10, 4)*{}="A"; (-7, 8)*{}="B"; (-4, 4)*{}="C";
(-10, -4)*{}="D"; (-7, -8)*{}="E"; (-4, -4)*{}="F";
(10, 4)*{}="A2"; (7, 8)*{}="B2"; (4, 4)*{}="C2";
(10, -4)*{}="D2"; (7, -8)*{}="E2"; (4, -4)*{}="F2";
"A"*{\bullet}; "B"*{\bullet}; "C"*{\bullet};
"D"*{\bullet}; "E"*{\bullet}; "F"*{\bullet};
"A2"*{\bullet}; "B2"*{\bullet}; "C2"*{\bullet};
"D2"*{\bullet}; "E2"*{\bullet}; "F2"*{\bullet};
"A"; "B"; **\dir{-}; "B"; "C"; **\dir{-}; "C"; "A"; **\dir{-};
"D"; "E"; **\dir{-}; "E"; "F"; **\dir{-}; "F"; "D"; **\dir{-};
"A2"; "B2"; **\dir{-}; "B2"; "C2"; **\dir{-}; "C2"; "A2"; **\dir{-};
"D2"; "E2"; **\dir{-}; "E2"; "F2"; **\dir{-}; "F2"; "D2"; **\dir{-};
(-12, 10)*{}; (-2, 10); **\dir{-};
(-12, -10)*{}; (-2, -10); **\dir{-};
(-12, 10)*{}; (-12, -10); **\dir{-};
(-2, 10)*{}; (-2, -10); **\dir{-};
(12, 10)*{}; (2, 10); **\dir{-};
(12, -10)*{}; (2, -10); **\dir{-};
(12, 10)*{}; (12, -10); **\dir{-};
(2, 10)*{}; (2, -10); **\dir{-};
\end{xy}
\qquad = \qquad
\begin{xy}
(-10, 4)*{}="A"; (-7, 8)*{}="B"; (-4, 4)*{}="C";
(-10, -4)*{}="D"; (-7, -8)*{}="E"; (-4, -4)*{}="F";
(10, 4)*{}="A2"; (7, 8)*{}="B2"; (4, 4)*{}="C2";
(10, -4)*{}="D2"; (7, -8)*{}="E2"; (4, -4)*{}="F2";
"A"*{\bullet}; "B"*{\bullet}; "C"*{\bullet};
"D"*{\bullet}; "E"*{\bullet}; "F"*{\bullet};
"A2"*{\bullet}; "B2"*{\bullet}; "C2"*{\bullet};
"D2"*{\bullet}; "E2"*{\bullet}; "F2"*{\bullet};
"A"; "B"; **\dir{-}; "B"; "C"; **\dir{-}; "C"; "A"; **\dir{-};
"D"; "E"; **\dir{-}; "E"; "F"; **\dir{-}; "F"; "D"; **\dir{-};
"A2"; "B2"; **\dir{-}; "B2"; "C2"; **\dir{-}; "C2"; "A2"; **\dir{-};
"D2"; "E2"; **\dir{-}; "E2"; "F2"; **\dir{-}; "F2"; "D2"; **\dir{-};
(-12, 10)*{}; (12, 10); **\dir{-};
(-12, 2)*{}; (12, 2); **\dir{-};
(-12, 10)*{}; (-12, 2); **\dir{-};
(12, 10)*{}; (12, 2); **\dir{-};
(-12, -10)*{}; (12, -10); **\dir{-};
(-12, -2)*{}; (12, -2); **\dir{-};
(-12, -10)*{}; (-12, -2); **\dir{-};
(12, -10)*{}; (12, -2); **\dir{-};
\end{xy}
\end{displaymath}
Clearly the proof of Proposition~\ref{prop:oddmerge} does not apply to this
relation.  We do not know if this relation lies in the span of the merging
relations.

\subsection{The space $W_L'$}

Let $W_L''$ (resp. $W_L'$) be the quotient of $\wt{W}_L$ by merging relations
(resp. merging and odd cycle exchange relations).  Consider the maps
\begin{displaymath}
\wt{W}_L \to W_L'' \to W_L' \to W_L.
\end{displaymath}
The first two maps are surjective; the third is surjective after 2 is inverted.
Proposition~\ref{prop:merge} is equivalent to the map $W_L'' \to W_L$ being an
isomorphism over $\Z[\half]$ when $L$ has cardinality at least 10 but not 12.
We have already shown that $W_L'' \to W_L'$ is an isomorphism over $\Z[\half]$
when the cardinality of $L$ is at least 10 but not 12.  Thus it suffices to
show that $W_L' \to W_L$ is an isomorphism over $\Z[\half]$ when the
cardinality of $L$ is at least 10 and not 12.  In fact, we will prove this
even when the cardinality of $L$ is equal to 12.

To prove this result, we need a way to think about elements of $W_L'$.  The
basic idea is that in quotienting $\wt{W}_L$ by the merging and odd cycle
exchange relations we have completely forgotten the partition data.  We now
elaborate on this idea.

\begin{definition}
(a) A degree two graph is \emph{forbidden} if it is connected or the union of
two odd cycles.  Graphs which are not forbidden are \emph{allowable}.
\newline
(b) Let $\Gamma$ be an allowable graph.  A pair of edges $(e, e')$ of $\Gamma$
is \emph{allowable} if it meets the following conditions:
\begin{itemize}
\item If $\Gamma$ is the union of two cycles then $e$ and $e'$ lie
in the same cycle.
\item If $\Gamma$ is the union of three cycles, two of which are
odd, then $e$ and $e'$ lie in cycles of the same parity.
\item Otherwise, there is no condition.
\end{itemize}
(c) A Pl\"ucker relation on an allowable pair of edges is called an
\emph{allowable Pl\"ucker relation}.
\end{definition}

Note that applying an allowable Pl\"ucker relation to an allowable graph
results in allowable graphs.  The description of $W_L'$ which we seek is the
following.  The proof is straightforward, and left to the reader.

\begin{proposition}
(a) Let $\Gamma$ be an allowable graph.  Choose a partition $\ms{U}$ of $L$
which is closed with respect to $\Gamma$.  Then the image of $(\Gamma, \ms{U})$
under the map $\wt{W}_L \to W_L'$ is independent of choice of $\ms{U}$.  We
denote the image by $X'_{\Gamma}$\\
(b) The $X'_{\Gamma}$ with $\Gamma$ allowable span $W'_L$ and satisfy the
allowable Pl\"ucker relations.\\
(c) The map $W'_L \to W_L$ takes $X'_{\Gamma}$ to $X_{\Gamma}$.
\end{proposition}

The statement that $W'_L \to W_L$ is an isomorphism over $\Z[\half]$ may now be
rephrased as the following statement:  given a collection of allowable graphs
$\{\Gamma_i\}$ such that $\sum a_i X_{\Gamma_i}=0$ when all Pl\"ucker relations
are allowed, the identity still holds if we restrict ourselves to the allowable
Pl\"ucker relations (and are allowed to invert 2).  We will prove this by
establishing a form of the straightening algorithm for the space $W_L'$.

\begin{remark}
In statement (b) of the proposition, it is in fact true that the allowable
Pl\"ucker relations generate all the linear relations between the
$X'_{\Gamma}$.  We will not make use of this fact, and  so do not provide a proof.
In fact, it will follow from our proof that $W'_L \to W_L$ is an isomorphism:
our proof shows that the map from the space of allowable graphs modulo
allowable Pl\"ucker relations to $W_L$ is an isomorphism.
\end{remark}

\subsection{Quasi-planar graphs}

Fix once and for all an embedding of $L$ into the unit circle in the plane.
We thus have a notation of planar graphs and Kempe's basis theorem
(Proposition~\ref{prop:straight}).  It is clear that allowable planar graphs
cannot form a basis for $W'_L$, since $W'_L \to W_L$ is surjective and there
are planar graphs which are forbidden.  Thus to prove something analogous to
Kempe's theorem for $W'_L$ we need to allow some non-planar graphs.  We choose
to allow graphs which are non-planar in the most mild manner possible:

\begin{definition}
We say that a graph $\Gamma$ on $L$ is \emph{quasi-planar} if either (1) it is
planar; or (2) $\Gamma$ has a doubled edge $e$ such that it becomes planar
when $e$ is removed, and $e$ crosses exactly two edges.
\end{definition}

We call any edge $e$ as in (2) a \emph{distinguished doubled edge} of $\Gamma$.
In any quasi-planar graph there are at most two distinguished doubled edges and
usually only one.  A distinguished doubled edge meets exactly one other
cycle of $\Gamma$.

\begin{figure}[!ht]
\begin{center}
\leavevmode
\begin{xy}
(0, 15)*{}="A"; (10.606, 10.606)*{}="B";
(15, 0)*{}="C"; (10.606, -10.606)*{}="D";
(0, -15)*{}="E"; (-10.606, -10.606)*{}="F";
(-15, 0)*{}="G"; (-10.606, 10.606)*{}="H";
(0, 0)*\xycircle(15, 15){};
"A"*{\bullet}; "B"*{\bullet}; "C"*{\bullet}; "D"*{\bullet};
"E"*{\bullet}; "F"*{\bullet}; "G"*{\bullet}; "H"*{\bullet};
"E"; "H"; **\dir{=};
"G"; "A"; **\dir{-}; "A"; "F"; **\dir{-}; "F"; "G"; **\dir{-};
"B"; "C"; **\dir{-}; "C"; "D"; **\dir{-}; "D"; "B"; **\dir{-};
\end{xy}
\hskip 2cm
\begin{xy}
(0, 15)*{}="A"; (10.606, 10.606)*{}="B";
(15, 0)*{}="C"; (10.606, -10.606)*{}="D";
(0, -15)*{}="E"; (-10.606, -10.606)*{}="F";
(-15, 0)*{}="G"; (-10.606, 10.606)*{}="H";
(0, 0)*\xycircle(15, 15){};
"A"*{\bullet}; "B"*{\bullet}; "C"*{\bullet}; "D"*{\bullet};
"E"*{\bullet}; "F"*{\bullet}; "G"*{\bullet}; "H"*{\bullet};
"A"; "E"; **\dir{-}; "E"; "F"; **\dir{-}; "F"; "G"; **\dir{-};
"G"; "H"; **\dir{-}; "H"; "A"; **\dir{-};
"B"; "C"; **\dir{-}; "C"; "D"; **\dir{-}; "D"; "B"; **\dir{-};
\end{xy}
\end{center}
\caption{A quasi-planar graph and its associated planar graph.
\label{f:quasi-planar}}
\end{figure}
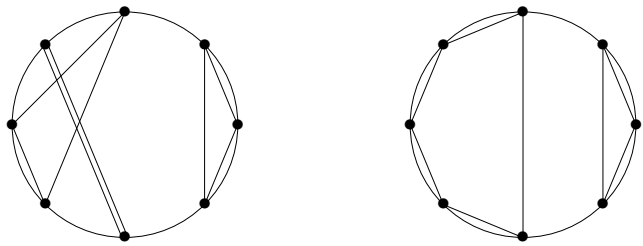

Let $\Gamma$ be a quasi-planar graph.  We can obtain a planar graph from
$\Gamma$ by replacing the distinguished doubled edge and the unique cycle it
crosses by the unique planar cycle on the same set of vertices.  We call this
graph the \emph{associated planar graph} to $\Gamma$.  See
Figure~\ref{f:quasi-planar} for an example.  We say that two quasi-planar
graphs are \emph{equivalent} if they have the same associated planar graph.

\begin{remark}
There is a sign ambiguity in this definition since we have not said how to
orient the edges in the associated planar graph.  This will not matter for
our purposes.
\end{remark}

\begin{lemma}
Assume the cardinality of $L$ is at least 8.  Then every planar graph is
associated to some allowable quasi-planar graph.
\end{lemma}

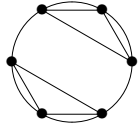
\begin{figure}[!ht]
\begin{displaymath}
\begin{xy}
(0, 8)*{}="A"; (6.928, 4)*{}="B"; (6.928, -4)*{}="C";
(0, -8)*{}="D"; (-6.928, -4)*{}="E"; (-6.928, 4)*{}="F";
(-4, -6.928)*{}="A"; (-8, 0)*{}="B"; (-4, 6.928)*{}="C";
(4, 6.928)*{}="D"; (8, 0)*{}="E"; (4, -6.928)*{}="F";
(0, 0)*\xycircle(8, 8){};
"A"*{\bullet}; "B"*{\bullet}; "C"*{\bullet};
"D"*{\bullet}; "E"*{\bullet}; "F"*{\bullet};
"A"; "B"; **\dir{-}; "B"; "F"; **\dir{-}; "F"; "A"; **\dir{-};
"C"; "D"; **\dir{-}; "D"; "E"; **\dir{-}; "E"; "C"; **\dir{-};
\end{xy}
\end{displaymath}
\caption{The unique planar graph which is not associated to any quasi-planar
graph.}
\label{f-no-qp}
\end{figure}

The proof of this lemma is straightforward.  See Figure~\ref{f-no-qp} for
a counterexample when $L$ has cardinality 6.

Let $\Gamma$ be a quasi-planar graph and let $n$ be the number of cycles in
$\Gamma$.  We define the \emph{level} of $\Gamma$ to be $n$ if $\Gamma$ is
planar and $n-1$ otherwise.  Thus in the quasi-planar case we just omit
the doubled edge in our count of cycles.  Note that $\Gamma$ and its associated
planar graph have the same level.

We define a decreasing filtration on $W_L$ (resp. $W_L'$) by letting $F^i W_L$
(resp.  $F^i W_L'$) be the span of all planar graphs (resp. allowable
quasi-planar graphs) with level at least $i$.  It is clear that both
filtrations are decreasing and separated.  The filtration $F^i W_L$ is also
exhaustive: $F^0 W_L=WL$ since the planar graphs span $W_L$.  It is not clear
that $F^0 W_L'=W_L'$, but we will prove this below.

For a planar graph $\Gamma$ of level $i$ we let $\ol{X}_\Gamma$ denote
the image of $X_{\Gamma}$ in $F^iW_L/F^{i+1}W_L$.  For an allowable
quasi-planar graph $\Gamma$ of level $i$ we let $\ol{X}'_{\Gamma}$ denote the
image of $X'_{\Gamma}$ in $F^iW_L'/F^{i+1}W_L'$.

\begin{proposition}
Let $\Gamma$ be an allowable quasi-planar graph of level $i$ and let
$\Gamma'$ be its associated planar graph, which also has level $i$.  Then
$X_{\Gamma}=\pm 2 X_{\Gamma'}$ modulo $F^{i-1} W_L$.  Thus the map
$W_L' \to W_L$ is filtration preserving and the image of $\ol{X}'_{\Gamma}$
under $\gr{W_L'} \to \gr{W_L}$ spans the same space as $\ol{X}_{\Gamma}$ if
2 is inverted.
\end{proposition}

\begin{proof}
The statement is obvious if $\Gamma$ is planar itself.  Thus assume that
$\Gamma$ is non-planar and let $e$ be a distinguished doubled edge.
We apply identity (I1) of the appendix to $\Gamma$, where vertices 1
and 4 belong to $e$ and the other two edges belong to the cycle which $e$
meets.  On the right side of the identity the first graph is simply $\Gamma'$.
The remaining graphs have more cycles, and thus belong to $F^{i+1} W_L$.
\end{proof}

As a consequence, we have:

\begin{corollary}
\label{cor:qp}
Let $L$ be an even set of cardinality at least 8.  Then over $\Z[\half]$:
\begin{itemize}
\item[(a)] Allowable quasi-planar graphs span $W_L$.
\item[(b)] Equivalent allowable quasi-planar graphs span the same subspace of
$\gr{W_L}$.
\end{itemize}
\end{corollary}

In the following two sections, we establish the following proposition:

\begin{proposition}
\label{prop:qp}
Let $L$ be an even set of cardinality at least 10.  Then over $\Z[\half]$:
\begin{itemize}
\item[(a)] Allowable quasi-planar graphs span $W_L'$.
\item[(b)] Equivalent allowable quasi-planar graphs span the same subspace
of $\gr{W_L'}$.
\end{itemize}
\end{proposition}

Proposition~\ref{prop:qp}, Corollary~\ref{cor:qp} and
Proposition~\ref{prop:straight} imply that $\gr{W_L'} \to \gr{W_L}$ is an
isomorphism.  Consequently, $W_L' \to W_L$ is an isomorphism.  Thus
Proposition~\ref{prop:merge} is implied by Proposition~\ref{prop:qp}.  (To see
this in another way, Proposition~\ref{prop:qp} shows that the dimension of
$W_L'$ at most the number of equivalence classes of allowable quasi-planar
graphs, which is equal to the number of planar graphs, which is equal to the
dimension of $W_L$ by Proposition~\ref{prop:straight}.  Thus $\dim{W_L'} \le
\dim{W_L}$.  Since $W_L' \to W_L$ is surjective it is therefore an
isomorphism.)

\section{Proof of Proposition~\ref{prop:qp}(\protect\textrm{a})}

\noindent
{\it We assume throughout this section that the cardinality of $L$ is even
and at least 10.  The word ``graph'' will always mean ``regular graph of
degree two.'' We also work over $\Z[\half]$ throughout this section.}
\vskip 1.2ex

We now prove Proposition~\ref{prop:qp}(a).  We will often use the
identities of the appendix.  It is worth noting that (I2) and (I3) hold
unconditionally, that is, all Pl\"ucker relations involved are always
allowable.  The reason for this is that these identities involve only a single
cycle.  We begin with a definition.

\begin{definition}
A vertex in a graph is \emph{special} if after deleting all edges containing
it the graph is planar.  A cycle in a graph is called \emph{special} if it
contains a special vertex.  A graph is \emph{semi-planar} if it contains
a special vertex.
\end{definition}

A quasi-planar graph is semi-planar.  A general semi-planar graph looks
like the following:
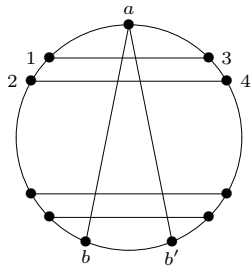
\begin{displaymath}
\begin{xy}
(0, 15)*{}="A"; (5.740, -13.858)*{}="B1"; (-5.740, -13.858)*{}="B2";
(-10.606, 10.606)*{}="C"; (10.606, 10.606)*{}="D";
(-12.99, 7.5)*{}="E"; (12.99, 7.5)*{}="F";
(-12.99, -7.5)*{}="G"; (12.99, -7.5)*{}="H";
(-10.606, -10.606)*{}="I"; (10.606, -10.606)*{}="J";
(0, 0)*\xycircle(15, 15){};
"A"*{\bullet}; "B1"*{\bullet}; "B2"*{\bullet};
"C"*{\bullet}; "D"*{\bullet}; "E"*{\bullet}; "F"*{\bullet};
"G"*{\bullet}; "H"*{\bullet}; "I"*{\bullet}; "J"*{\bullet};
"A"; "B1"; **\dir{-};
"A"; "B2"; **\dir{-};
"C"; "D"; **\dir{-};
"E"; "F"; **\dir{-};
"G"; "H"; **\dir{-};
"I"; "J"; **\dir{-};
(-13, 10.606)*{\ss 1}; (-15.5, 7.5)*{\ss 2};
(13, 10.606)*{\ss 3}; (15.5, 7.5)*{\ss 4};
(0, 17)*{\ss a};
(-5.740, -15.858)*{\ss b};
(5.740, -15.858)*{\ss b'};
\end{xy}
\end{displaymath}
Here $a$ is a special vertex.  The vertices 1, 2, 3 and 4 are part of a planar
cycle, as are the four vertices at the bottom of the circle.  There may be
other planar cycles which either do or do not intersect one of the two edges
containing $a$.

\begin{lemma}
The $X'_{\Gamma}$ with $\Gamma$ an allowable semi-planar graph span
$W'_L$ over $\Z[\half]$.
\end{lemma}

\begin{proof}
Elements of the form $X'_{\Gamma}$, with $\Gamma$ an allowable graph
containing a doubled edge, span $W'_L$ over $\Z[\half]$.  This can be seen
using the results
of \cite[\S 6.3]{hmsv2} or the identity (I2).  Now, if $\Gamma$ has a doubled
edge $e$ then one can hold $e$ fixed and apply the straightening algorithm
to the rest of $\Gamma$.  All the resulting graphs have the property that
they become planar when $e$ is removed.  They are thus semi-planar.  Hence
$X'_{\Gamma}$ has been expressed as a sum of semi-planar graphs.
\end{proof}

Let $\Gamma$ be a semi-planar graph and let $a$ be a special vertex.  The
vertex $a$ meets two edges.  If one starts at $a$ and follows one of these
edges one will intersect some of the cycles of $\Gamma$.  In fact, the set
of cycles met, and the order in which they are met, is independent of the
edge chosen.  We call these cycles the \emph{skewered cycles}.  We call
the last skewered cycle encountered (starting from $a$) \emph{extreme}.
We now make a few simple observations.

\begin{lemma}
\label{lem:aqp}
Let $\Gamma$ be an allowable graph.  Assume one of the following conditions
is satisfied:
\begin{itemize}
\item[(a)]  $\Gamma$ has a doubled edge $e$ which crosses at most two other
edges.
\item[(b)] $\Gamma$ has four distinct vertices $a$, $b$, $c$, and $d$ such
that $\uedge{ab}$, $\uedge{bc}$ and $\uedge{cd}$ are edges in $\Gamma$ which
do not cross any edge.
\item[(c)] $\Gamma$ contains an $n$-cycle ($n \ge 4$) with the property that
no edge in the $n$-cycle crosses any edge.
\item[(d)] $\Gamma$ contains two 3-cycles with the property that no edge in
either cycle crosses any edge.
\item[(e)] $\Gamma$ is semi-planar and has a special $n$-cycle with $n \ge 4$.
\item[(f)] $\Gamma$ is semi-planar and has a skewered $n$-cycle with 
$n \ne 3, 4$.
\end{itemize}
Then $X'_{\Gamma}$ is a sum of allowable quasi-planar graphs over $\Z[\half]$.
\end{lemma}

\begin{proof}
(a) Fix the edge $e$ and apply the straightening algorithm to the remainder
of $\Gamma$.  The key point is that in the straightening algorithm one only
has to Pl\"ucker edges which cross and therefore no new crossings with $e$ are
introduced.  Thus in all the graphs resulting from the straightening algorithm
$e$ will still only cross two edges.  These graphs are therefore quasi-planar.

(b) First consider the case where $\uedge{ad}$ is an edge in $\Gamma$.  It is
clear that this edge is not crossed by any other edge.  Holding the 4-cycle
$abcd$ fixed and applying the straightening algorithm to the rest of $\Gamma$
write $X'_{\Gamma}$ as a sum of allowable planar graphs.  Now consider the case
where $\uedge{ad}$ is not an edge in $\Gamma$.  Then $a$ is connected
to a unique vertex $x$ other than $b$ and $d$ is connected to a unique
vertex $y$ other than $c$.  It may be that $x=y$.  We now apply identity
(I3) with 1, 2, 3, 4, 5, 6 taken to be $x$, $a$, $b$, $c$, $d$, $y$.  The
first graph on the right side can be handled by the first case.  The
remaining graphs can be handled by (a) of this lemma.

(c) If $n$ is even simply hold the $n$-cycle fixed and apply the straightening
algorithm to the rest of the graph.  If $n$ is odd then it is $\ge 5$.  Apply
identity (I2) to five consecutive vertices in the cycle.  The first graph
has an even cycle which does not cross any edge and has already been handled.
The remaining graphs are handled by (a) of this lemma.

(d) Simply hold the two 3-cycles fixed and apply the straightening algorithm
to the remainder of the graph.  The key point is that since we have fixed
\emph{two} 3-cycles all the Pl\"ucker relations in the straightening algorithm
are allowable.

(e) If there is a special cycle of length at least five then this cycle
has four consecutive vertices none of which is the special vertex.  It
can be written as a sum of allowable quasi-planar graphs by part (b) of this
lemma.  Now consider the case where there is a special cycle of length 4.
Apply the identity:
\begin{equation}
\label{eq:sqr}
2 \hskip 1ex
\begin{xy}
(0, 8)*{}="A"; (8, 0)*{}="B"; (0, -8)*{}="C"; (-8, 0)*{}="D";
(0, 10)*{\ss a}; (10, 0)*{\ss b}; (0, -10)*{\ss c}; (-10, 0)*{\ss d};
(0, 0)*\xycircle(8, 8){};
"A"*{\bullet}; "B"*{\bullet}; "C"*{\bullet}; "D"*{\bullet};
"A"; "B"; **\dir{-};
"B"; "C"; **\dir{-};
"C"; "D"; **\dir{-};
"D"; "A"; **\dir{-};
\end{xy}
\quad = \quad
\begin{xy}
(0, 8)*{}="A"; (8, 0)*{}="B"; (0, -8)*{}="C"; (-8, 0)*{}="D";
(0, 10)*{\ss a}; (10, 0)*{\ss b}; (0, -10)*{\ss c}; (-10, 0)*{\ss d};
(0, 0)*\xycircle(8, 8){};
"A"*{\bullet}; "B"*{\bullet}; "C"*{\bullet}; "D"*{\bullet};
"A"; "B"; **\dir{=};
"C"; "D"; **\dir{=};
\end{xy}
\quad + \quad
\begin{xy}
(0, 8)*{}="A"; (8, 0)*{}="B"; (0, -8)*{}="C"; (-8, 0)*{}="D";
(0, 10)*{\ss a}; (10, 0)*{\ss b}; (0, -10)*{\ss c}; (-10, 0)*{\ss d};
(0, 0)*\xycircle(8, 8){};
"A"*{\bullet}; "B"*{\bullet}; "C"*{\bullet}; "D"*{\bullet};
"A"; "C"; **\dir{=};
"B"; "D"; **\dir{=};
\end{xy}
\quad + \quad
\begin{xy}
(0, 8)*{}="A"; (8, 0)*{}="B"; (0, -8)*{}="C"; (-8, 0)*{}="D";
(0, 10)*{\ss a}; (10, 0)*{\ss b}; (0, -10)*{\ss c}; (-10, 0)*{\ss d};
(0, 0)*\xycircle(8, 8){};
"A"*{\bullet}; "B"*{\bullet}; "C"*{\bullet}; "D"*{\bullet};
"A"; "D"; **\dir{=};
"B"; "C"; **\dir{=};
\end{xy}
\end{equation}
This identity can be gotten by applying the straightening algorithm to the
middle term on the right.  Here we take $a$ to be the special vertex.  Thus
$\uedge{bc}$ and $\uedge{cd}$ are not crossed by any edge.  Each of the
graphs on the right can thus be handled by (a) of this lemma; precisely we
take the edge $e$ in (a) to be $\uedge{cd}$, $\uedge{bd}$ and $\uedge{bc}$.

(f) Consider a skewered $n$-cycle.  If $n=2$ we are done by (a).  Thus
assume $n \ge 5$.  The special cycle splits the $n$ vertices
of this cycle into two sets.  If one of these sets has cardinality at least
four then the result follows from (b) of this lemma.  This is always the case
when $n \ge 7$.  Say now $n=6$.  We must handle the case where the special
cycle splits the 6 vertices into two sets of size 3.  Apply identity (I2)
where 2, 3, 4 all lie to one side of the special cycle.  The first graph
on the right has two 3-cycles of the sort that can be handled by (d) of
this lemma.  The remaining three graphs can be handled by (a) of this lemma.
Now say $n=5$.  We must handle the case where the special cycle splits
the 5 vertices into a set of size 3 and a set of size 2.  Apply (I2) where
again 2, 3 and 4 lie to one side of the special cycle.  The first graph now
has a doubled edge and can be handled by (a); the remaining three graphs
can be handled with (a) as before.
\end{proof}

The following lemma completes the proof of Proposition~\ref{prop:qp}(a).

\begin{lemma}
Let $\Gamma$ be an allowable semi-planar graph.  Then $X'_{\Gamma}$ is a
sum of allowable quasi-planar graphs over $\Z[\half]$.
\end{lemma}

\begin{proof}
The idea is to proceed inductively on the number of skewered cycles.  Thus
let $\Gamma$ be given.  In the generic case, simple apply the straightening
algorithm to the special cycle and the extreme skewered edge.  All the
resulting graphs will be semi-planar and have fewer skewered cycles.  This
procedure does not work in some small cases because the straightening
algorithm may involve forbidden Pl\"ucker relations.  In fact, this procedure
fails only in the following two cases:  $\Gamma$ has two even cycles; or
$\Gamma$ has two odd cycles, one even cycle and the special cycle and extreme
skewered cycle have opposite parity.  We now handle each of these cases,
breaking the second into two subcases.

First assume that $\Gamma$ has only two even cycles.  By
Lemma~\ref{lem:aqp}(e), if the special cycle has length at least four then we
are done.  Thus we may assume that the special cycle is a doubled edge.  But
now the other cycle has length at least 8 (since there are at least 10
vertices by assumption) and the result follows from Lemma~\ref{lem:aqp}(f).

Now assume that $\Gamma$ has three cycles, that the special cycle has
even length and the extreme skewered cycle has odd length.  By
Lemma~\ref{lem:aqp} it suffices to treat the case where the special
cycle is a doubled edge and the extreme skewered cycle is a 3-cycle.  Thus
the remaining cycle has length at least 5.  If it is skewered we are done
by Lemma~\ref{lem:aqp}(f) otherwise we are done by Lemma~\ref{lem:aqp}(c).

Finally assume that $\Gamma$ has three cycles, that the special cycle has
odd length and the extreme skewered cycle has even length.  By
Lemma~\ref{lem:aqp} if suffices to treat the case where the special cycle
is a 3-cycle and the extreme skewered cycle is a 4-cycle.  If there are
at least 12 vertices in total then there is an $n$-cycle with $n \ge 5$ and
we are done as before.  However, if there are 10 vertices then we have to
do some work.  The special cycle and extreme skewered cycle look like:
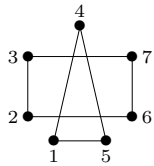
\begin{displaymath}
\begin{xy}
(0, 8)*{}="A"; (6.928, 4)*{}="B"; (6.928, -4)*{}="C";
(3.471, -7.208)*{}="D"; (-3.471, -7.208)*{}="E";
(-6.928, -4)*{}="F"; (-6.928, 4)*{}="G";
(-3.471, -9.208)*{\ss 1};
(-8.928, -4)*{\ss 2};
(-8.928, 4)*{\ss 3};
(0, 10)*{\ss 4};
(3.471, -9.208)*{\ss 5};
(8.928, -4)*{\ss 6};
(8.928, 4)*{\ss 7};
"A"*{\bullet}; "B"*{\bullet}; "C"*{\bullet};
"D"*{\bullet}; "E"*{\bullet}; "F"*{\bullet}; "G"*{\bullet};
"A"; "D"; **\dir{-};
"A"; "E"; **\dir{-};
"D"; "E"; **\dir{-};
"B"; "G"; **\dir{-};
"B"; "C"; **\dir{-};
"C"; "F"; **\dir{-};
"F"; "G"; **\dir{-};
\end{xy}
\end{displaymath}
Here 4 is the special vertex.
We cannot apply the straightening algorithm, as this would involve forbidden
Pl\"ucker relations.  Instead, we apply \eqref{eq:sqr} to the square above.
We now have
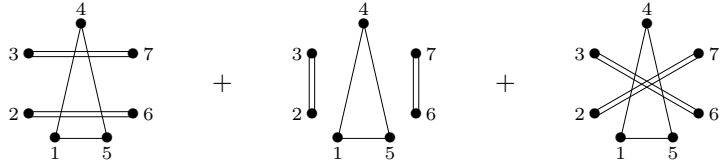
\begin{displaymath}
\begin{xy}
(0, 8)*{}="A"; (6.928, 4)*{}="B"; (6.928, -4)*{}="C";
(3.471, -7.208)*{}="D"; (-3.471, -7.208)*{}="E";
(-6.928, -4)*{}="F"; (-6.928, 4)*{}="G";
(-3.471, -9.208)*{\ss 1};
(-8.928, -4)*{\ss 2};
(-8.928, 4)*{\ss 3};
(0, 10)*{\ss 4};
(3.471, -9.208)*{\ss 5};
(8.928, -4)*{\ss 6};
(8.928, 4)*{\ss 7};
"A"*{\bullet}; "B"*{\bullet}; "C"*{\bullet};
"D"*{\bullet}; "E"*{\bullet}; "F"*{\bullet}; "G"*{\bullet};
"A"; "D"; **\dir{-};
"A"; "E"; **\dir{-};
"D"; "E"; **\dir{-};
"B"; "G"; **\dir{=};
"C"; "F"; **\dir{=};
\end{xy}
\qquad + \qquad
\begin{xy}
(0, 8)*{}="A"; (6.928, 4)*{}="B"; (6.928, -4)*{}="C";
(3.471, -7.208)*{}="D"; (-3.471, -7.208)*{}="E";
(-6.928, -4)*{}="F"; (-6.928, 4)*{}="G";
(-3.471, -9.208)*{\ss 1};
(-8.928, -4)*{\ss 2};
(-8.928, 4)*{\ss 3};
(0, 10)*{\ss 4};
(3.471, -9.208)*{\ss 5};
(8.928, -4)*{\ss 6};
(8.928, 4)*{\ss 7};
"A"*{\bullet}; "B"*{\bullet}; "C"*{\bullet};
"D"*{\bullet}; "E"*{\bullet}; "F"*{\bullet}; "G"*{\bullet};
"A"; "D"; **\dir{-};
"A"; "E"; **\dir{-};
"D"; "E"; **\dir{-};
"B"; "C"; **\dir{=};
"F"; "G"; **\dir{=};
\end{xy}
\qquad + \qquad
\begin{xy}
(0, 8)*{}="A"; (6.928, 4)*{}="B"; (6.928, -4)*{}="C";
(3.471, -7.208)*{}="D"; (-3.471, -7.208)*{}="E";
(-6.928, -4)*{}="F"; (-6.928, 4)*{}="G";
(-3.471, -9.208)*{\ss 1};
(-8.928, -4)*{\ss 2};
(-8.928, 4)*{\ss 3};
(0, 10)*{\ss 4};
(3.471, -9.208)*{\ss 5};
(8.928, -4)*{\ss 6};
(8.928, 4)*{\ss 7};
"A"*{\bullet}; "B"*{\bullet}; "C"*{\bullet};
"D"*{\bullet}; "E"*{\bullet}; "F"*{\bullet}; "G"*{\bullet};
"A"; "D"; **\dir{-};
"A"; "E"; **\dir{-};
"D"; "E"; **\dir{-};
"B"; "F"; **\dir{=};
"C"; "G"; **\dir{=};
\end{xy}
\end{displaymath}
The first two graphs can be handled by Lemma~\ref{lem:aqp}(a).  In the
third graph, we hold $\uedge{27}$ fixed and apply the straightening
algorithm to the remaining edges.  This, of course, involves only allowable
Pl\"ucker relations.  The result is:
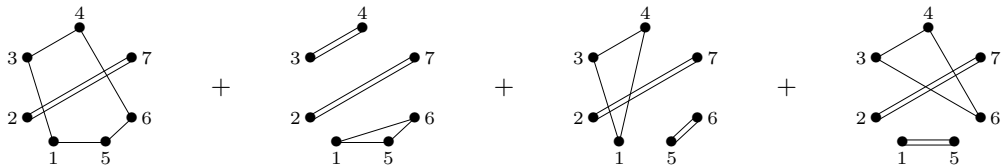
\begin{displaymath}
\begin{xy}
(0, 8)*{}="A"; (6.928, 4)*{}="B"; (6.928, -4)*{}="C";
(3.471, -7.208)*{}="D"; (-3.471, -7.208)*{}="E";
(-6.928, -4)*{}="F"; (-6.928, 4)*{}="G";
(-3.471, -9.208)*{\ss 1};
(-8.928, -4)*{\ss 2};
(-8.928, 4)*{\ss 3};
(0, 10)*{\ss 4};
(3.471, -9.208)*{\ss 5};
(8.928, -4)*{\ss 6};
(8.928, 4)*{\ss 7};
"A"*{\bullet}; "B"*{\bullet}; "C"*{\bullet};
"D"*{\bullet}; "E"*{\bullet}; "F"*{\bullet}; "G"*{\bullet};
"B"; "F"; **\dir{=};
"A"; "C"; **\dir{-};
"C"; "D"; **\dir{-};
"D"; "E"; **\dir{-};
"E"; "G"; **\dir{-};
"G"; "A"; **\dir{-};
\end{xy}
\qquad + \qquad
\begin{xy}
(0, 8)*{}="A"; (6.928, 4)*{}="B"; (6.928, -4)*{}="C";
(3.471, -7.208)*{}="D"; (-3.471, -7.208)*{}="E";
(-6.928, -4)*{}="F"; (-6.928, 4)*{}="G";
(-3.471, -9.208)*{\ss 1};
(-8.928, -4)*{\ss 2};
(-8.928, 4)*{\ss 3};
(0, 10)*{\ss 4};
(3.471, -9.208)*{\ss 5};
(8.928, -4)*{\ss 6};
(8.928, 4)*{\ss 7};
"A"*{\bullet}; "B"*{\bullet}; "C"*{\bullet};
"D"*{\bullet}; "E"*{\bullet}; "F"*{\bullet}; "G"*{\bullet};
"B"; "F"; **\dir{=};
"A"; "G"; **\dir{=};
"C"; "D"; **\dir{-};
"D"; "E"; **\dir{-};
"E"; "C"; **\dir{-};
\end{xy}
\qquad + \qquad
\begin{xy}
(0, 8)*{}="A"; (6.928, 4)*{}="B"; (6.928, -4)*{}="C";
(3.471, -7.208)*{}="D"; (-3.471, -7.208)*{}="E";
(-6.928, -4)*{}="F"; (-6.928, 4)*{}="G";
(-3.471, -9.208)*{\ss 1};
(-8.928, -4)*{\ss 2};
(-8.928, 4)*{\ss 3};
(0, 10)*{\ss 4};
(3.471, -9.208)*{\ss 5};
(8.928, -4)*{\ss 6};
(8.928, 4)*{\ss 7};
"A"*{\bullet}; "B"*{\bullet}; "C"*{\bullet};
"D"*{\bullet}; "E"*{\bullet}; "F"*{\bullet}; "G"*{\bullet};
"B"; "F"; **\dir{=};
"C"; "D"; **\dir{=};
"E"; "G"; **\dir{-};
"G"; "A"; **\dir{-};
"A"; "E"; **\dir{-};
\end{xy}
\qquad + \qquad
\begin{xy}
(0, 8)*{}="A"; (6.928, 4)*{}="B"; (6.928, -4)*{}="C";
(3.471, -7.208)*{}="D"; (-3.471, -7.208)*{}="E";
(-6.928, -4)*{}="F"; (-6.928, 4)*{}="G";
(-3.471, -9.208)*{\ss 1};
(-8.928, -4)*{\ss 2};
(-8.928, 4)*{\ss 3};
(0, 10)*{\ss 4};
(3.471, -9.208)*{\ss 5};
(8.928, -4)*{\ss 6};
(8.928, 4)*{\ss 7};
"A"*{\bullet}; "B"*{\bullet}; "C"*{\bullet};
"D"*{\bullet}; "E"*{\bullet}; "F"*{\bullet}; "G"*{\bullet};
"B"; "F"; **\dir{=};
"D"; "E"; **\dir{=};
"A"; "C"; **\dir{-};
"C"; "G"; **\dir{-};
"G"; "A"; **\dir{-};
\end{xy}
\end{displaymath}
Each of these graphs has 4 as a special vertex.  Thus in all of them
the doubled edge $\uedge{27}$ only crosses the edges which are drawn.
Each of the above graphs is therefore a sum of allowable quasi-planar
graphs by Lemma~\ref{lem:aqp}(a).
\end{proof}

\section{Proof of Proposition~\ref{prop:qp}(b)}

\noindent
{\it We maintain the assumptions from the beginning of the last section.}
\vskip 1.2ex

We now prove Proposition~\ref{prop:qp}(b).  We handle the three cases of level
1, level 2 and level at least 3 separately.  We begin by considering graphs of
level 3 or more.

\begin{lemma}
\label{lem-lev3}
Let $\Gamma$ be an allowable quasi-planar graph of level at least 3 and
let $\Gamma'$ be its associated planar graph.  Then $\Gamma'$ is allowable
and $\ol{X}'_{\Gamma}$ and $\ol{X}'_{\Gamma'}$ span the
same subspace of $\gr{W_L'}$.  In particular, any two equivalent allowable
quasi-planar graphs of level at least 3 span the same subspace of
$\gr{W_L'}$.
\end{lemma}

\begin{proof}
Simply apply identity (I1) to the doubled edge in $\Gamma$.  The first graph is
$\Gamma'$ and the rest have higher level.
\end{proof}

We now we consider graphs of level 2.  There are three possibilities for
such a graph $\Gamma$:
\begin{itemize}
\item $\Gamma$ has two odd cycles and one doubled edge which crosses one of
the odd cycles.
\item $\Gamma$ has two even cycles and one doubled edge which crosses one
of the even cycles.
\item $\Gamma$ has two even cycles.
\end{itemize}
In the latter two cases, the same reasoning as Lemma~\ref{lem-lev3}
applies.  We are thus reduced to the first case.  We call such a graph
``odd.''

Let $\Gamma$ be an odd level 2 allowable quasi-planar graph.  Let $e$ be the
2-cycle, let $C$ be the cycle which crosses $e$ and let $C'$ be the other
cycle.  The vertices used in $C$ and $e$ form an interval, as do the vertices
used in $C'$.  We call these two intervals the \emph{distinguished intervals}
of $\Gamma$.  It follows from the definition that two odd graphs are equivalent
if and only if they have the same distinguished intervals.

\begin{lemma}
Let $\Gamma$ be an odd level 2 allowable quasi-planar graph and let $c$, $d$
and $e$ be consecutively connected vertices which do not form a 3-cycle.  Let
$\Gamma'$ be the unique odd graph which is equivalent to $\Gamma$ and for which
$\uedge{ce}$ is the doubled edge.  Then $\ol{X}'_{\Gamma}$ and $\ol{X}'_{\Gamma'}$
span the same subspace of $\gr{W_L'}$.
\end{lemma}

\begin{proof}
We apply (I2) to $\Gamma$ with 2, 3 and 4 taken to be $c$, $d$ and $e$.  The
last two graphs on the right of the identity are allowable planar since they
have 2-cycles on consecutive vertices.  The first graph can also be written as
a sum of allowable planar graphs (it is quasi-planar and has four cycles).
This leaves the second graph.  This graph has two doubled edges:  the original
one and $\uedge{ce}$.  We apply identity (I1) to the original doubled edge.  The
first graph on the right is $\Gamma'$ and the remaining have higher level.
This gives an expression $\ol{X}'_{\Gamma}=\pm 2 \ol{X}'_{\Gamma'}$ and
proves the lemma.
\end{proof}

\begin{corollary}
Two equivalent allowable quasi-planar graphs of level 2 span the same subspace
of $\gr{W_L'}$.
\end{corollary}

\begin{proof}
Let $\Gamma$ be a given odd level 2 graph and let $\ol{\Gamma}$ be its
associated planar graph.  Let $n<m$ be the lengths of the two
cycles in $\ol{\Gamma}$; these are odd numbers whose sum is at least 10.  We
first consider the case where $n$ and $m$ are each at least 5.  Pick three
consecutively connected vertices $c_0$, $d_0$ and $e_0$ in $\ol{\Gamma}$.
Now, let $\uedge{ce}$ be the distinguished doubled edge of $\Gamma$.  If $c$ and
$e$ lie in the same component of $\ol{\Gamma}$ as $c_0$ apply the lemma to
move the doubled edge to the other component.  Now one can again apply the
lemma to move the doubled edge of $\Gamma$ to $\uedge{c_0 d_0}$.  This shows that
any allowable quasi-planar graphs with associated planar $\ol{\Gamma}$ spans
the same space in $\gr{W'_L}$ as the unique allowable quasi-planar graph
associated to $\ol{\Gamma}$ and with distinguished doubled edge $\uedge{c_0 d_0}$.

We now consider the case where $n=3$ so that $m \ge 7$.  The distinguished
doubled edge must lie in the $m$-cycle of $\Gamma$.  Pick three consecutively
connected vertices $c_0$, $d_0$ and $e_0$ in $\ol{\Gamma}$ which lie in the
$m$-cycle.  Using the lemma we can first move the doubled edge of $\Gamma$
far away from $c_0$ and $e_0$ and then move it to be exactly $\uedge{c_0 d_0}$.
By the same reasoning at the end of the previous paragraph this completes
the proof.
\end{proof}

We finish by considering allowable quasi-planar graphs of level 1.  These are
graphs which have a single doubled edge and exactly one other cycle, which is
planar and which meets the doubled edge.  Such a graph is determined
by which two vertices belong to the doubled edge.

\begin{lemma}
\label{lem1}
Let $\Gamma$ be an allowable quasi-planar  graph of level 1 and let $c$, $d$
and $e$ be consecutively connected vertices.  Let $\Gamma'$ be the unique
allowable quasi-planar level 1 graph which has $\uedge{c e}$ for its doubled edge.
Then $\ol{X}'_{\Gamma}$ and $\ol{X}'_{\Gamma'}$ span the same subspace of
$\gr{W_L'}$.
\end{lemma}

\begin{proof}
We apply (I2) to $\Gamma$ with 2, 3, 4 taken to be $c$, $d$ and $e$.  Note
that $c$, $d$ and $e$ necessarily lie in an $n$-cycle with $n \ge 8$, so (I2)
can indeed be applied.  The first graph on the right side has level at least
two.  The third and fourth graphs have a 2-cycle on
consecutive vertices and so can be written as a sum of allowable planars.
Now, the second graph has $\uedge{ce}$ as a doubled edge and two other cycles:
the original doubled edge and remainder of the original long cycle.  We apply
identity (I1) to the original doubled edge and the two edges it crosses.
The first graph on the right is $\Gamma'$ and the rest have higher level.
We thus obtain an expression $\ol{X}'_{\Gamma}=\pm 2\ol{X}'_{\Gamma'}$, which
establishes the lemma.
\end{proof}

\begin{corollary}
Two equivalent allowable quasi-planar graphs of level 1 span the same
subspace of $\gr{W_L'}$.
\end{corollary}

\begin{proof}
This goes like the second paragraph of the proof of the previous corollary.
\end{proof}

    %%%%%%%%%%%%%%
    %% APPENDIX %%
    %%%%%%%%%%%%%%

\appendix
\section{Long identities}
\label{app}

In all the following identities we use the convention that edges point from
smaller numbers to larger numbers. 

\vskip 3ex

\par\noindent
Identity (I1):
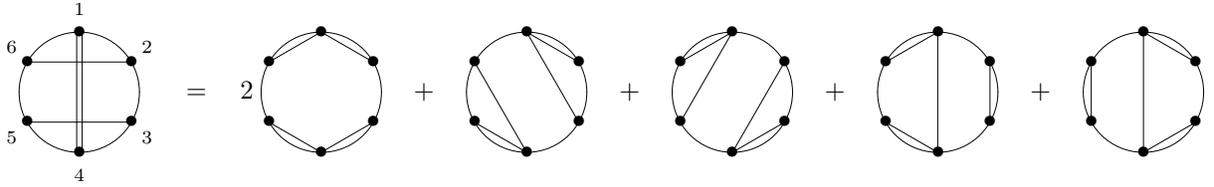
\begin{displaymath}
\begin{xy}
(0, 8)*{}="A"; (6.928, 4)*{}="B"; (6.928, -4)*{}="C";
(0, -8)*{}="D"; (-6.928, -4)*{}="E"; (-6.928, 4)*{}="F";
(0, 0)*\xycircle(8, 8){};
"A"*{\bullet}; "B"*{\bullet}; "C"*{\bullet};
"D"*{\bullet}; "E"*{\bullet}; "F"*{\bullet};
(0, 11)*{\ss 1}; (9, 6)*{\ss 2}; (9, -6)*{\ss 3};
(0, -11)*{\ss 4}; (-9, -6)*{\ss 5}; (-9, 6)*{\ss 6};
"A"; "D"; **\dir{=};
"B"; "F"; **\dir{-}; "C"; "E"; **\dir{-};
\end{xy}
\quad = \quad
2 \;
\begin{xy}
(0, 8)*{}="A"; (6.928, 4)*{}="B"; (6.928, -4)*{}="C";
(0, -8)*{}="D"; (-6.928, -4)*{}="E"; (-6.928, 4)*{}="F";
(0, 0)*\xycircle(8, 8){};
"A"*{\bullet}; "B"*{\bullet}; "C"*{\bullet};
"D"*{\bullet}; "E"*{\bullet}; "F"*{\bullet};
"A"; "B"; **\dir{-}; "C"; "D"; **\dir{-};
"D"; "E"; **\dir{-}; "F"; "A"; **\dir{-};
\end{xy}
\quad + \quad
\begin{xy}
(0, 8)*{}="A"; (6.928, 4)*{}="B"; (6.928, -4)*{}="C";
(0, -8)*{}="D"; (-6.928, -4)*{}="E"; (-6.928, 4)*{}="F";
(0, 0)*\xycircle(8, 8){};
"A"*{\bullet}; "B"*{\bullet}; "C"*{\bullet};
"D"*{\bullet}; "E"*{\bullet}; "F"*{\bullet};
"A"; "B"; **\dir{-}; "A"; "C"; **\dir{-};
"D"; "E"; **\dir{-}; "D"; "F"; **\dir{-};
\end{xy}
\quad + \quad
\begin{xy}
(0, 8)*{}="A"; (6.928, 4)*{}="B"; (6.928, -4)*{}="C";
(0, -8)*{}="D"; (-6.928, -4)*{}="E"; (-6.928, 4)*{}="F";
(0, 0)*\xycircle(8, 8){};
"A"*{\bullet}; "B"*{\bullet}; "C"*{\bullet};
"D"*{\bullet}; "E"*{\bullet}; "F"*{\bullet};
"A"; "E"; **\dir{-}; "A"; "F"; **\dir{-};
"D"; "B"; **\dir{-}; "D"; "C"; **\dir{-};
\end{xy}
\quad + \quad
\begin{xy}
(0, 8)*{}="A"; (6.928, 4)*{}="B"; (6.928, -4)*{}="C";
(0, -8)*{}="D"; (-6.928, -4)*{}="E"; (-6.928, 4)*{}="F";
(0, 0)*\xycircle(8, 8){};
"A"*{\bullet}; "B"*{\bullet}; "C"*{\bullet};
"D"*{\bullet}; "E"*{\bullet}; "F"*{\bullet};
"A"; "D"; **\dir{-}; "B"; "C"; **\dir{-};
"A"; "F"; **\dir{-}; "D"; "E"; **\dir{-};
\end{xy}
\quad + \quad
\begin{xy}
(0, 8)*{}="A"; (6.928, 4)*{}="B"; (6.928, -4)*{}="C";
(0, -8)*{}="D"; (-6.928, -4)*{}="E"; (-6.928, 4)*{}="F";
(0, 0)*\xycircle(8, 8){};
"A"*{\bullet}; "B"*{\bullet}; "C"*{\bullet};
"D"*{\bullet}; "E"*{\bullet}; "F"*{\bullet};
"A"; "D"; **\dir{-}; "E"; "F"; **\dir{-};
"A"; "B"; **\dir{-}; "D"; "C"; **\dir{-};
\end{xy}
\end{displaymath}
\vskip 1ex
\par\noindent
This identity can be proved by applying the straightening algorithm to
the left side.

\vskip 3ex

\par\noindent
Identity (I2):
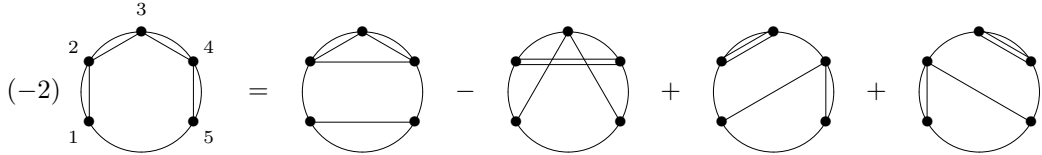
\begin{displaymath}
(-2) \;
\begin{xy}
(-6.928, -4)*{}="A"; (-6.928, 4)*{}="B"; (0, 8)*{}="C";
(6.928, 4)*{}="D"; (6.928, -4)*{}="E";
(0, 0)*\xycircle(8, 8){};
"A"*{\bullet}; "B"*{\bullet}; "C"*{\bullet};
"D"*{\bullet}; "E"*{\bullet};
(-9, -6)*{\ss 1}; (-9, 6)*{\ss 2}; (0, 11)*{\ss 3};
(9, 6)*{\ss 4}; (9, -6)*{\ss 5};
"A"; "B"; **\dir{-}; "B"; "C"; **\dir{-};
"C"; "D"; **\dir{-}; "D"; "E"; **\dir{-};
\end{xy}
\quad = \quad
\begin{xy}
(-6.928, -4)*{}="A"; (-6.928, 4)*{}="B"; (0, 8)*{}="C";
(6.928, 4)*{}="D"; (6.928, -4)*{}="E";
(0, 0)*\xycircle(8, 8){};
"A"*{\bullet}; "B"*{\bullet}; "C"*{\bullet};
"D"*{\bullet}; "E"*{\bullet};
"A"; "E"; **\dir{-}; "B"; "C"; **\dir{-};
"C"; "D"; **\dir{-}; "D"; "B"; **\dir{-};
\end{xy}
\quad - \quad
\begin{xy}
(-6.928, -4)*{}="A"; (-6.928, 4)*{}="B"; (0, 8)*{}="C";
(6.928, 4)*{}="D"; (6.928, -4)*{}="E";
(0, 0)*\xycircle(8, 8){};
"A"*{\bullet}; "B"*{\bullet}; "C"*{\bullet};
"D"*{\bullet}; "E"*{\bullet};
"A"; "C"; **\dir{-}; "C"; "E"; **\dir{-};
"B"; "D"; **\dir{=};
\end{xy}
\quad + \quad
\begin{xy}
(-6.928, -4)*{}="A"; (-6.928, 4)*{}="B"; (0, 8)*{}="C";
(6.928, 4)*{}="D"; (6.928, -4)*{}="E";
(0, 0)*\xycircle(8, 8){};
"A"*{\bullet}; "B"*{\bullet}; "C"*{\bullet};
"D"*{\bullet}; "E"*{\bullet};
"A"; "D"; **\dir{-}; "D"; "E"; **\dir{-};
"B"; "C"; **\dir{=};
\end{xy}
\quad + \quad
\begin{xy}
(-6.928, -4)*{}="A"; (-6.928, 4)*{}="B"; (0, 8)*{}="C";
(6.928, 4)*{}="D"; (6.928, -4)*{}="E";
(0, 0)*\xycircle(8, 8){};
"A"*{\bullet}; "B"*{\bullet}; "C"*{\bullet};
"D"*{\bullet}; "E"*{\bullet};
"A"; "B"; **\dir{-}; "B"; "E"; **\dir{-};
"C"; "D"; **\dir{=};
\end{xy}
\end{displaymath}
\vskip 1ex
\par\noindent
This identity can be proved by applying the straightening algorithm to
the second graph on the right side

\vskip 3ex

\par\noindent
Identity (I3):
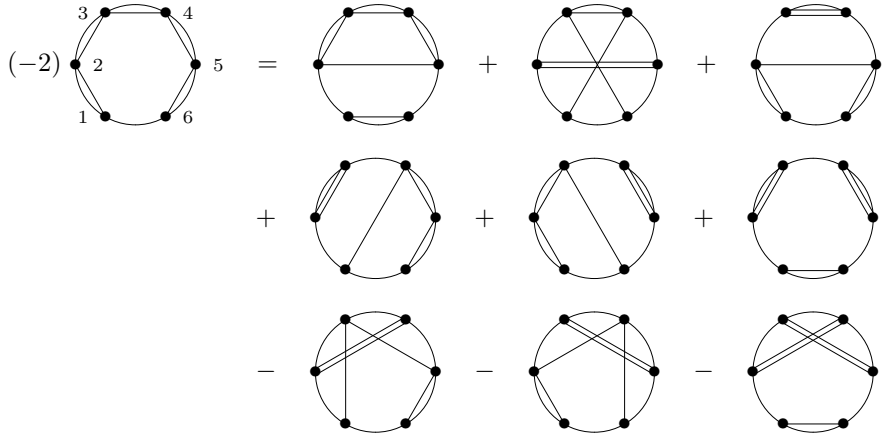
\begin{displaymath}
\begin{split}
(-2) \;
\begin{xy}
(-4, -6.928)*{}="A"; (-8, 0)*{}="B"; (-4, 6.928)*{}="C";
(4, 6.928)*{}="D"; (8, 0)*{}="E"; (4, -6.928)*{}="F";
(0, 0)*\xycircle(8, 8){};
"A"*{\bullet}; "B"*{\bullet}; "C"*{\bullet};
"D"*{\bullet}; "E"*{\bullet}; "F"*{\bullet};
(-7, -6.928)*{\ss 1}; (-5, 0)*{\ss 2}; (-7, 6.928)*{\ss 3};
(7, 6.928)*{\ss 4}; (11, 0)*{\ss 5}; (7, -6.928)*{\ss 6};
"A"; "B"; **\dir{-}; "B"; "C"; **\dir{-};
"C"; "D"; **\dir{-}; "D"; "E"; **\dir{-};
"E"; "F"; **\dir{-};
\end{xy}
\quad &= \quad
\begin{xy}
(-4, -6.928)*{}="A"; (-8, 0)*{}="B"; (-4, 6.928)*{}="C";
(4, 6.928)*{}="D"; (8, 0)*{}="E"; (4, -6.928)*{}="F";
(0, 0)*\xycircle(8, 8){};
"A"*{\bullet}; "B"*{\bullet}; "C"*{\bullet};
"D"*{\bullet}; "E"*{\bullet}; "F"*{\bullet};
"A"; "F"; **\dir{-}; "B"; "C"; **\dir{-};
"C"; "D"; **\dir{-}; "D"; "E"; **\dir{-};
"E"; "B"; **\dir{-};
\end{xy}
\quad + \quad
\begin{xy}
(-4, -6.928)*{}="A"; (-8, 0)*{}="B"; (-4, 6.928)*{}="C";
(4, 6.928)*{}="D"; (8, 0)*{}="E"; (4, -6.928)*{}="F";
(0, 0)*\xycircle(8, 8){};
"A"*{\bullet}; "B"*{\bullet}; "C"*{\bullet};
"D"*{\bullet}; "E"*{\bullet}; "F"*{\bullet};
"A"; "D"; **\dir{-}; "D"; "C"; **\dir{-};
"C"; "F"; **\dir{-}; "B"; "E"; **\dir{=};
\end{xy}
\quad + \quad
\begin{xy}
(-4, -6.928)*{}="A"; (-8, 0)*{}="B"; (-4, 6.928)*{}="C";
(4, 6.928)*{}="D"; (8, 0)*{}="E"; (4, -6.928)*{}="F";
(0, 0)*\xycircle(8, 8){};
"A"*{\bullet}; "B"*{\bullet}; "C"*{\bullet};
"D"*{\bullet}; "E"*{\bullet}; "F"*{\bullet};
"A"; "B"; **\dir{-}; "B"; "E"; **\dir{-};
"E"; "F"; **\dir{-}; "C"; "D"; **\dir{=};
\end{xy}
\\[2ex] &+ \quad
\begin{xy}
(-4, -6.928)*{}="A"; (-8, 0)*{}="B"; (-4, 6.928)*{}="C";
(4, 6.928)*{}="D"; (8, 0)*{}="E"; (4, -6.928)*{}="F";
(0, 0)*\xycircle(8, 8){};
"A"*{\bullet}; "B"*{\bullet}; "C"*{\bullet};
"D"*{\bullet}; "E"*{\bullet}; "F"*{\bullet};
"A"; "D"; **\dir{-}; "D"; "E"; **\dir{-};
"E"; "F"; **\dir{-}; "B"; "C"; **\dir{=};
\end{xy}
\quad + \quad
\begin{xy}
(-4, -6.928)*{}="A"; (-8, 0)*{}="B"; (-4, 6.928)*{}="C";
(4, 6.928)*{}="D"; (8, 0)*{}="E"; (4, -6.928)*{}="F";
(0, 0)*\xycircle(8, 8){};
"A"*{\bullet}; "B"*{\bullet}; "C"*{\bullet};
"D"*{\bullet}; "E"*{\bullet}; "F"*{\bullet};
"A"; "B"; **\dir{-}; "B"; "C"; **\dir{-};
"C"; "F"; **\dir{-}; "D"; "E"; **\dir{=};
\end{xy}
\quad + \quad
\begin{xy}
(-4, -6.928)*{}="A"; (-8, 0)*{}="B"; (-4, 6.928)*{}="C";
(4, 6.928)*{}="D"; (8, 0)*{}="E"; (4, -6.928)*{}="F";
(0, 0)*\xycircle(8, 8){};
"A"*{\bullet}; "B"*{\bullet}; "C"*{\bullet};
"D"*{\bullet}; "E"*{\bullet}; "F"*{\bullet};
"A"; "F"; **\dir{-};
"B"; "C"; **\dir{=}; "D"; "E"; **\dir{=};
\end{xy}
\\[2ex] &- \quad
\begin{xy}
(-4, -6.928)*{}="A"; (-8, 0)*{}="B"; (-4, 6.928)*{}="C";
(4, 6.928)*{}="D"; (8, 0)*{}="E"; (4, -6.928)*{}="F";
(0, 0)*\xycircle(8, 8){};
"A"*{\bullet}; "B"*{\bullet}; "C"*{\bullet};
"D"*{\bullet}; "E"*{\bullet}; "F"*{\bullet};
"A"; "C"; **\dir{-}; "C"; "E"; **\dir{-};
"E"; "F"; **\dir{-}; "B"; "D"; **\dir{=};
\end{xy}
\quad - \quad
\begin{xy}
(-4, -6.928)*{}="A"; (-8, 0)*{}="B"; (-4, 6.928)*{}="C";
(4, 6.928)*{}="D"; (8, 0)*{}="E"; (4, -6.928)*{}="F";
(0, 0)*\xycircle(8, 8){};
"A"*{\bullet}; "B"*{\bullet}; "C"*{\bullet};
"D"*{\bullet}; "E"*{\bullet}; "F"*{\bullet};
"A"; "B"; **\dir{-}; "B"; "D"; **\dir{-};
"D"; "F"; **\dir{-}; "C"; "E"; **\dir{=};
\end{xy}
\quad - \quad
\begin{xy}
(-4, -6.928)*{}="A"; (-8, 0)*{}="B"; (-4, 6.928)*{}="C";
(4, 6.928)*{}="D"; (8, 0)*{}="E"; (4, -6.928)*{}="F";
(0, 0)*\xycircle(8, 8){};
"A"*{\bullet}; "B"*{\bullet}; "C"*{\bullet};
"D"*{\bullet}; "E"*{\bullet}; "F"*{\bullet};
"A"; "F"; **\dir{-};
"B"; "D"; **\dir{=}; "C"; "E"; **\dir{=};
\end{xy}
\end{split}
\end{displaymath}
\vskip 1ex
\par\noindent
This identity is proved in \cite[\S 6.3]{hmsv2}.  It can be proved by applying
the straightening algorithm to the three graphs on the last line.

\vskip 3ex

    %%%%%%%%%%%%%%%%
    %% REFERENCES %%
    %%%%%%%%%%%%%%%%

\end{document}